\documentclass[11pt]{amsart}\usepackage{amscd,amsmath,amsfonts,amssymb}

\newcommand{\la}{\lambda}

\newcommand{\Ll}{\mathcal{L}}

\newcommand{\CR}{$\rm{CR}$}

\newcommand{\CC}{\mathbb{C}}
\newcommand{\HH}{\mathbb{H}}
\newcommand{\RR}{\mathbb{R}}

\numberwithin{equation}{section}

\newtheorem{pr}{Proposition}[section]

\theoremstyle{definition}

\newtheorem{re}{Remark}[section]

\begin{document}

\title{\CR-submanifolds. A class of examples.}

\author[Liviu Ornea]{Liviu Ornea}

\thanks{I am pleased to announce financial support from the University of New Mexico at Albuquerque during the preparation of this note.\\[.5cm]
}

 \address{\newline
\noindent University of Bucharest,\newline Faculty of
Mathematics,\newline 14 Academiei str.,\newline  70109 Bucharest,
Romania}

\email{ Liviu.Ornea@imar.ro, lornea@gta.math.unibuc.ro}

\begin{abstract}
I present a class of examples of \CR-submanifolds of manifolds endowed with different structures, obtained as level sets of momentum maps associated to specific Hamiltonian actions.
\end{abstract}
\maketitle

\begin{center}
\emph{Dedicated to Tudor Zamfirescu at his 60th birthday.}
\end{center}

\section{Motivation}
Tudor always attacks and can solve very concrete, geometrical problems, ones that you are able to see. I cannot. He doesn't like structures on manifolds. I do. So, what kind of mathematical gift could I possibly offer him? I thought that the less ``offensive'' would be some examples.

\bigskip

\CR-submanifolds were introduced first in K\"ahler geometry, as a bridge between complex and totally real submanifolds. In Bejancu's definition, roughly speaking, their tangent bundle splits into a complex part of constant dimension and a totally real part, orthogonal to the first one. After its introduction, in \cite{be1}, the definition was soon extended to other ambient spaces and gave rise to a large amount of literature. The monographs \cite{yk1}, \cite{yk2},  \cite{be2} and others gather only part of the (most significant) results, since their publication many others appeared, the list remains open.

What structures presumably admit such submanifolds? All is needed is (1) a (pseudo)-Riemannian metric and (2) a compatible 1--1 tensor field with respect to which the invariant and antiinvariant distributions (generalizations of complex and totally real ones) should be defined. In general, this tensor field is chosen to satisfy a certain polynomial equation. Of course, for an arbitrary tensor, one loses the motivating property that a \CR-submanifold of an Hermitian manifold is \CR in the sense of complex variables, nevertheless their geometry might be interesting.

All the same, there are not so many natural examples of \CR-submanifolds. The aim of this very short note is to make explicit  a general class of examples, some of them already mentioned by passing in previous papers like \cite{gop}, \cite{dror}, examples naturally arising whenever one has a manifold structured as above and a group action giving rise to a momentum map.

\section{Structured manifolds and Hamiltonian actions}

\subsection{Definitions and results} 
Let $(M,g,F)$ be a (semi-)Riemannian manifold endowed with a 1--1 tensor field $F$ compatible with $g$ in the following sense:
$$g(X,FY)+g(FX,Y)=0.$$
In other words, I require that 
$$\omega(X,Y):= g(X,FY)$$
be a $2$-form. I shall suppose $\omega$ of maximal rank.

I call such a triple a \emph{structured manifold}. Note that no integrability of $F$ is required. Obvious particular cases are:
\begin{itemize}
\item almost Hermitian manifolds ($F^2=-1$);
\item almost contact manifolds ($F^2=-1+\eta\otimes\xi$, where $\eta$ and $\xi$ are dual with respect to $g$ and $F\xi=0$);
\item $f$-manifolds ($F^3+F=0$);
\item almost product manifolds ($F^2=1$ and $g$ has non-trivial index).
\end{itemize}  

A submanifold $N$ of a structured manifold $M$ is called a \emph{\CR-submanifold} if its tangent bundle splits into two orthogonal distributions $D$ and $D^\perp$ such that:
\begin{enumerate}
\item $TN=D\oplus D^\perp$.
\item $FD=D$ and $FD^\perp\subseteq T^\perp N$.
\end{enumerate}

Traditionally, the name \CR-submanifold is reserved for a Hermitian ambient, contact-\CR\ or semi-invariant are the corresponding names in almost contact manifolds etc.

Suppose $G$ is a (connected) group of automorphisms of a  structured manifold as above with Lie algebra $\mathfrak{g}$. This means that $G$ acts by isometries and the flow of any fundamental field induced by the action commutes with $F$: $\Ll_Xg=0$ and $\Ll_XF=0$ for any fundamental field $X$, where $\Ll_X$ denotes the Lie derivative on the direction $X$. I shall say that such a group contains automorphisms of the structure.

Suppose, moreover, that this action has an associated equivariant momentum map. This means the existence of a map $\mu:M\rightarrow \mathfrak{g}^*$ satisfying the following two conditions:
\begin{enumerate}
\item $d\mu (\zeta)=\zeta_M \rfloor \omega$, where $\zeta\in \mathfrak{g}$ and $\zeta_M$ is the induced fundamental field.
\item $\mu(a\cdot x)=a\cdot \mu(x)$, where the left- (resp. right-) hand side $\cdot$ represents the action of $G$ on $M$ (resp. the coadjoint action on $\mathfrak{g}^*$).
\end{enumerate}

If $\omega$ is closed, the existence of the momentum map is equivalent with the action being (strongly) Hamiltonian. But there are cases in which a momentum map exists even if the 2-form is not closed: this happens, for example,  on locally conformal K\"ahler manifolds (cf. \cite{gop}), on  contact and Sasakian manifolds (cf. \cite{dror}). Even more radical a situation appears when one drops up even the existence of $g$, like in Joyce's constructions of hypercomplex quotient (\cite{joy}).

If a momentum map exists, one can prove, by just slightly generalizing the symplectic and K\"ahler situation, the following

\begin{pr}
Let $(M,g,F)$ be a structured manifold, $G$ a connected group of automorphisms of the structure. If an associated momentum map exists for which $0\in\mathfrak{g}$ is a regular value, then $\mu^{-1}(0)$ is a \CR-submanifold of $(M,g,F)$.
\end{pr}

\begin{proof}
Note first that $X\in T_xM$ is tangent to $\mu^{-1}(0)$ if and only if $d\mu(X)=0$, hence, by the definition of the momentum map, if and only if $\omega(X,\cdot)=0$.

If $G\cdot x$ represents the orbit of $G$ through $x$, then it easy to see that there exists the orthogonal splitting
$$T_x\mu^{-1}(0)=T_x(G\cdot x)\oplus H_x$$
with $H_x$ $F$-invariant: indeed this follows from the equivalence $X$ from $T_xM$ stays in $H_x$ if and only if $g(X,Y)=0$ and $\omega(X,Y)=g(X,FY)=0$ for all $Y\in T_x(G\cdot x)$. All in all we have $g(FX,Y)=0$ and $\omega(FX,Y)=0$, which proves the assertion.

Finally, if $X\in  T_x(G\cdot x)$ and $Y\in T_x\mu^{-1}(0)$ arbitrary, then $g(Y,FX)=\omega(Y,X)=d\mu(Y)(X)=0$, hence $T_x(G\cdot x)$ is a totally real distribution.

One then takes $D_x=H_x$ and $D^\perp=T_x(G\cdot x)$.   
\end{proof}

\begin{re}\label{nonzero}
The same proof does not work, in general, for non-zero regular values $\zeta$ of $\mu$. But it will work when $G_\zeta=G$, in particular for abelian $G$ (cf. \cite{br} for the K\"ahler case). Another situation when one can work with non-zero regular values is the one described in \cite{dror}: there, for $\zeta\neq 0$ in $\mathfrak{g}^*$, one looks at  $\mu^{-1}(\RR_+\zeta)$ which proves to be a \CR-submanifold.
\end{re}

\begin{re}
The well-known result stating that the anti-invariant distribution $D^\perp$ is integrable in l.c.K. (in particular, K\"ahler) manifolds, cf. \cite{blch}, extended afterwards to other classes of ambient spaces, is here a mere consequence of the fact that this distribution is isomorphic with the Lie algebra $\mathfrak{g}$: for this class of examples, the $D^\perp$ distribution is integrable, no matter the integrability properties of the ambient structure. Note that there are examples, cf. \emph{loc.cit}, where the totally real distribution is not integrable.
\end{re}

\begin{re}
When the topological conditions in the specific reduction theorems are satisfied, then the reduced space $M_0$ inherits the structure of $M$ and the natural projection $\mu^{-1}(0)\rightarrow M_0$ becomes a $\rm{CR}$-submersion, see \cite{kob} for the K\"ahler case, \cite{dror} for the Sasakian one.
\end{re}  

\begin{re}
The ${\rm CR}$-submanifolds of this type are \emph{generic} in the sense that $T_x^\perp\mu^{-1}(0)=FD_x^\perp$ and also in Chen's acception: the $F$-invariant distribution $D$ is the maximal one with this property. But if one is able to consider non-zero reduction, as in \cite{dror} (see Remark \ref{nonzero}), the resulted \CR-submanifolds are no longer generic.
\end{re}

\begin{re} An open and challenging problem is to describe properties of these \CR-sub\-ma\-ni\-folds (as they are usually studied in submanifolds theory) in terms of properties of the action. If the action satisfies the good topological assumptions in order that the quotient exist, this amounts to describing (Riemannian) 
properties of the quotient in terms of group action. But this setting is more general, the quotient might not exist, not even as an orbifold; in this case one would have to deal with distributions and transverse geometry.
\end{re}  

\begin{re}
The result also shows a less obvious relation between  \CR-submanifolds and symmetries of a structured manifold. In principle, the larger the automorphism group and the larger the number of its subgroups, more chances for the existence of \CR-submanifolds. On the other hand, clearly there exist \CR-submanifolds which are not level sets of momentum maps.
\end{re}

\subsubsection{Action by conformalities}
A more general situation (from the viewpoint of the group action) is the following. Let $(M,F,g)$ be a structured manifold such that $\omega$ is locally conformal with a closed form (this is equivalent with the existence of a closed 1-form $\theta$ such that $d\omega=\theta\wedge\omega$). For the moment, mainly are  studied the locally conformal K\"ahler (l.c.K.) and locally conformally cosymplectic (l.c.c.) manifolds and their submanifolds. Note that on the universal cover of such a manifold, the pull-back metric $\tilde g$ is globally conformal with a metric $h$ whith closed fundamental form. 

One can slightly modify the above procedure to obtain examples of \CR-submanifolds. The change consists in allowing the group $G$ act by conformalities (not mere isometries) and define an appropriate "twisted Hamiltonian" action (this amounts to replacing the operator of exterior derivative $d$ with $d^\theta:=d+\theta\wedge$ and accordingly modifying all definitions; see \cite{gop} and \cite{gopp} for details in the l.c.K. case). As above, one may prove:
\begin{pr}
Let $(M,F,g)$ be a structured manifold such that $\omega$ is locally conformal with a closed form and let $G$ be a connected group acting by conformalities preserving $F$. If the action is twisted Hamiltonian with momentum map $\mu$, then $\mu^{-1}(0)$ is a \CR-submanifold of  $(M,F,g)$.
\end{pr}
 
\subsection{Compatibility among various types of \CR-submanifolds}
Much studied in the literature is the following situation: Let $\tilde M$ be a Sasakian manifold with regular Reeb field $\xi$ and $M=M/\xi$ its K\"ahler base, $\pi$ the canonical projection (a principal $S^1$ fibration). If $\tilde N$ is a submanifold of $\tilde M$ which fibers in circles over a submanifold $N$ of $M$,  then $\tilde N$ is a contact \CR-submanifold if and only if $N$ is a \CR-submanifold. Various other properties of the two can now be related.

This picture naturally appears also in our setting. Let $\pi:\tilde M\rightarrow M$ be a submersion between structured manifolds, commuting with the respective structure tensors  (such a map is sometimes called $(\tilde F, F)$-holomorphic). Suppose $G$ is a group of automorphisms of both structures, commuting with $\pi$ and  providing respective momentum maps $\tilde \mu$ and $\mu$. We then have $\tilde\mu=\mu\circ\pi$ and it can be easily shown that $\pi$ restricts to a submersion between the \CR-submanifolds $\tilde\mu^{-1}(0)$ and $\mu^{-1}(0)$. 

\subsubsection{Action by conformalities}
Let now $(M,g,F)$ be a structured manifold with fundamental form locally conformal with a closed one. We let $G$ be a (connected) group of conformalities preserving $F$ acting on $M$. Let $\tilde M$ be a cover of $M$ on which $\theta$ is closed (in particular, $\tilde M$ can be the universal cover, but it is not necessary). Let also $\tilde G$ be a lift of $G$ to $\tilde M$ (see \cite{gopp} for a thorough discussion of the existence and properties of $\tilde G$) and let $\tilde G_0$ be its identity component.

The following result was proven in \cite{gop} for l.c.K. manifolds, but it  can be extended immediately to other structures:

\begin{pr} 
In the above setting, if the action of $G$ on $M$ is twisted Hamiltonian, then the action of $\tilde G_0$ on $\tilde M$ is Hamiltonian with respect to the metric $h$. In particular, a twisted Hamiltonian action by conformalities lifts to a Hamiltonian action by isometries.
\end{pr}

Hence, starting with a twisted Hamiltonian action on $(M,g,F)$, with momentum map $\mu$, one obtains a Hamiltonian action on $(\tilde M, h, F)$ with momentum map $\tilde \mu$. One then has a \CR-submanifold $\mu^{-1}(0)$ in $M$ a \CR-submanifold $\tilde \mu^{-1}(0)$ in $\tilde M$ and it is easily seen that:
\begin{pr}
The covering $\tilde M\rightarrow M$ restricts to a covering $\tilde\mu^{-1}(0)\rightarrow \mu^{-1}(0)$ of \CR-submanifolds.
\end{pr}

\section{Some examples}
\noindent {\bf 1.}~ Any real hypersurface of a structured manifold is a (trivial) example of \CR-submanifold, with 1-dimensional distribution $D^\perp$. 

{\bf 1.1.}~In \cite{gror} examples of this kind in Sasakian manifolds are obtained using circle actions.  For example, acting on $S^{2n-1}$ with its standard Sasakian structure as follows:
$$(e^{it}, z)\mapsto (e^{it\la_0}z_0,\ldots,e^{\la_nit}z_{n-1}),$$
with integer weights $\la_a$, the momentum map is $\mu(z)=\sum \la_a|z_a|^2$. To have a non-empty level set, one assumes the weights do not have all the same sign. In general, it is not an easy task to determine (at least up to homeomorphism) the level set. If $\la_0=\cdots =\la_k=a>0$ and $\la_{k+1}=\cdots =\la_{n-1}=-b<0$, then it is easy to identify $\mu^{-1}(0)$ with the product $S^{2k+1}(\sqrt{a/a+b})\times S^{2(n-k)-1}(\sqrt{b/a+b})$ which becomes a \CR-submanifold of $S^{2n-1}$.

{\bf 1.2.}~ Similar computations lead in \cite{gop} to \CR-submanifolds of l.c.K. manifolds, specifically of Hopf manifolds.

{\bf 1.3.}~ An example where the ambient space is para-Hermitian (metric almost product) appears in \cite{kond}. On $\RR^n\times \RR^n$ with coordinates $(x_j,y_j)$ one considers the almost product structure given by $F=\sum(\partial y_j\otimes dx_j+\partial x_j\otimes dy_j)$ and the pseudo-Riemannian metric $\sum (dx_j^2-dy_j^2)$. The associated 2-form is closed. Now one acts with $(\RR,+)$ by $(t, (x,y))\mapsto (x\cosh t+y\sinh t, x\sinh t+y\cosh t)$. The momentum map is $\mu(x,y)=\sum (x_j^2-y_j^2)$ and the hypersurface $\mu^{-1}(\frac 12)\cong S^{n-1,n}$ is a \CR-submanifold.

\noindent {\bf 2.}~ A non-trivial example of \CR-submanifold in the standard Sasakian sphere $S^{4n-1}$, associated to the $\mathrm{SU}(2)$-action given by right multiplication appears in \cite{gror'}. In quaternionic coordinates, the momentum map $\mu:S^{4n-1}\rightarrow \mathfrak{su}(2)^*=\RR^3$ is $\mu(q)=\sum q_ai\bar q_a$. The zero level set can be seen to be the homogeneous manifold
$\displaystyle \frac{\mathrm{SU}(n+1)}{\mathrm{SU}(n-1)}$, which is thus a \CR-submanifold of the sphere.

\noindent {\bf 3.}~ Another non-trivial example produced by a weighted action of $T^2$ on $S^7$ appears in \cite{dror}. Let $S^7$ be endowed with the standard Sasakian structure and act on it with $T^2$ as follows:
$$((e^{it_0},e^{it_1}),z)\mapsto
(e^{it_{0}\la_0}z_0,e^{it_{1}\la_1}z_1, z_2, z_3),$$
where $\la_0$, $\la_1$ are some real weights. The momentum map, written as a linear map on $\RR^2$ is:
$$\mu(z)=\langle(\la_0\vert z_0\vert^2,\la_1\vert
 z_1\vert^2),\cdot\rangle$$
 and its zero level set is $S^3$. 

\noindent {\bf 4.} According to Remark \ref{nonzero}, in the Sasakian case one is also able to perform non-zero reduction. Several examples of this type are given in \cite{dror} using actions of $T^2$ on spheres. For example, using different actions of $T^2$ on $S^7$ one explicitly obtains the following (non-generic) \CR-submanifolds: 
\begin{itemize}
\item $S^3$. Here the action is $((e^{it_0},e^{it_1}),z)\mapsto
(e^{it_{0}}z_0,e^{it_{0}}z_1, e^{it_{1}}z_2, e^{it_{1}}z_3)$, the momentum map is $\mu(z)=\langle (|z_0|^2+|z_1|^2, z_2|^2+|z_3|^2), \cdot\rangle$ and $\zeta=\langle(1,0),\cdot\rangle$. This example can be generalized to obtain \CR-submanifolds of $S^{2n-1}$.
\item $S^5$. Here the action is $((e^{it_0},e^{it_1}),z)\mapsto
(e^{it_{0}}z_0,e^{it_{1}}z_1, e^{it_{1}}z_2, e^{it_{1}}z_3)$, the momentum map reads $mu(z)=\langle (|z_0|^2+|z_1|^2, z_2|^2+|z_3|^2), \cdot\rangle$ and $\zeta=\langle(0, v), \cdot\rangle)$.
\item $S^1(\frac{1}{\sqrt{2}})\times (S^5({\sqrt{2})}\setminus S^1(\frac{1}{\sqrt{2}}))$, an open submanifold. Here the action is $((e^{it_0},e^{it_1}),z)\mapsto
(e^{-it_{0}}z_0,e^{it_{0}}z_1, e^{it_{1}}z_2, e^{it_{1}}z_3)$, the momentum map is $\mu(z)=\langle (|z_0|^2-|z_1|^2, z_2|^2+|z_3|^2), \cdot\rangle$ and $\zeta=\langle(1,1),\cdot\rangle$.
\item $S^7\cap (\CC^*\times A)$, where $A$ is the ellipsoid of equation $\mid z_1\mid^2(1+\la_1/\la_0)+\mid z_2\mid^2+\mid z_3\mid^2=1$, the $\la's$ being the weights of the action, see above.
\end{itemize}

\begin{re} In all the examples concerning Sasakian ambient spaces, the Reeb field is tangent to the constructed \CR-submanifold. Moreover, all the considered actions commute with the Reeb field (which, for the standard Sasakian structure of the sphere is regular) and hence, according to Section 2.2, the quotient by the Reeb flow of the \CR-submanifolds of the form $\mu^{-1}(0)$ are \CR-submanifolds of the complex projective space. 
\end{re} 

\noindent{\bf 5.} Let $T^2$ act on $\CC^4$ by
$$((s,t), z)\mapsto (s^mtz_1, tz_2, sz_3, sz_4).$$
One easily verifies that the action is K\"ahlerian with respect to the standard complex structure and the Fubini-Study metric (in fact, this is the action that produces as symplectic quotient the Hirzebruch surfaces, cf. \cite{aud}). The momentum map is
$$\mu(z)=\frac 12(m|z_1|^2+|z_3|^2+|z_4|^2, |z_1|^2+|z_2|^2), \quad |s|=|t|=1,$$
and one may verify that $(r,1)\in \RR^2=(\mathfrak{t}^2)^*$ is a regular value for $\mu$ provided $0\neq r\neq m$. Hence, for $r>m$,  $\mu^{-1}(r,1)$ is a $6$-dimensional \CR-submanifold of $\CC^4$ (cf. Remark \ref{nonzero}).

\noindent{\bf 6.} Let $(\CC^{n+k})^k$ be the space of complex matrices of type $(n+k, k)$ endowed with the symplectic form $\omega(X,Y)=\mathrm{Im}\mathrm{Tr}(\bar X^tY)$. The group $\mathrm{U}(k)$ acts naturally on this space with momentum map $\mu (A)=\bar A^tA$ (identifying $\mathfrak{u}(k)^*$ with the space of Hermitian matrices). The level set $\mu^{-1}(Id)$ is identified with the Stiefel manifold $V_k(\CC^{n+k})$ (cf. \cite{aud}) which is thus a \CR-submanifold of $(\CC^{n+k})^k$.

\noindent{\bf 7.} Other examples can be extracted from hyperk\"ahler reductions, fixing one of the K\"ahler structures of the ambient manifold. Here is one of them (cf. \cite{grg}), appearing in the generalization of the Taub-Nut metric. Let $M=\HH^n\times \HH^n$ with (quaternionic) coordinates $(q_a,w_a)$. We consider it as a K\"ahler manifold with respect to the flat product metric and the complex structure given by the left multiplication with $i$. Let $G=\RR^n=\{(t_j)\}$ acting on $M$ as follows:
$$q_a\mapsto q_ae^{it_a}, \quad w_a\mapsto w_a-\la_a^bt_b,$$
where $(\la_a^b)$ is a real non-degenerate $(m,m)$ matrix. The action is K\"ahlerian (in fact, hyperk\"ahlerian) and the momentum map with respect to the chosen complex structure is the $i$ component of the following map:
$$\mu_a(q_a,w_a)=\frac 12\sum q_ai\bar q_a+\frac 12\sum \la_a^b(w_b-\bar w_b), \quad a=1,\ldots,n.$$
 Note that the conjugation is the quaternionic one. The zero level set is then a \CR-submanifold.

\end{document}